\documentclass[12pt]{article}

\usepackage{mathrsfs}

\usepackage{amsmath,amsthm,amssymb}
\usepackage{graphicx}
\usepackage{vector}

\font\tencmmib=cmmib10 \skewchar\tencmmib '60
\newfam\cmmibfam
\textfont\cmmibfam=\tencmmib

\def\bbox{\quad\hbox{\vrule \vbox{\hrule \vskip2pt \hbox{\hskip2pt
\vbox{\hsize=1pt}\hskip2pt} \vskip2pt\hrule}\vrule}}
\def\lessim{\ \lower4pt\hbox{$
\buildrel{\displaystyle <}\over\sim$}\ }
\def\gessim{\ \lower4pt\hbox{$\buildrel{\displaystyle >}
\over\sim$}\ }

\def\C{{\cal C}}

\def\A{{\cal A}}

\def\I{{\cal I}}

\def\la{{\Bigl\langle}}
\def\ra{{\Bigr\rangle}}

\def\qed{\hfill\break\rightline{$\bbox$}}

\newcommand{\e}{\mathbb{E}}

\newcommand{\Reals}{\mathbb{R}}
\newcommand{\Natural}{\mathbb{N}}

\newcommand{\vsi}{{\vec{\sigma}}}
\newcommand{\vrho}{{\vec{\rho}}}

\newcommand{\s}{\boldsymbol{\sigma}}

\newtheorem{lemma}{\bf Lemma}

\newtheorem{theorem}{\bf Theorem}

\newtheorem{proposition}{\bf Proposition}

\newenvironment{Proof of lemma}{\noindent{\bf Proof of Lemma}}{\hfill$\Box$\newline}
\newenvironment{Proof of theorem}{\noindent{\bf Proof of Theorem}}{\hfill$\Box$\newline}
\newenvironment{Proof of theorems}{\noindent{\bf Proof of Theorems}}{\hfill$\Box$\newline}
\newenvironment{Proof of proposition}{\noindent{\bf Proof of Proposition}}{\hfill$\Box$\newline}
\newenvironment{Proof of propositions}{\noindent{\bf Proof of Propositions}}{\hfill$\Box$\newline}
\newenvironment{Proof of exercise}{\noindent{\it Proof of Exercise:}}{\hfill$\Box$}

\font\tencmmib=cmmib10 \skewchar\tencmmib '60
\newfam\cmmibfam
\textfont\cmmibfam=\tencmmib

\def\bbox{\quad\hbox{\vrule \vbox{\hrule \vskip2pt \hbox{\hskip2pt
\vbox{\hsize=1pt}\hskip2pt} \vskip2pt\hrule}\vrule}}
\def\lessim{\ \lower4pt\hbox{$
\buildrel{\displaystyle <}\over\sim$}\ }
\def\gessim{\ \lower4pt\hbox{$\buildrel{\displaystyle >}
\over\sim$}\ }

\def\go0{\to 0}

\def\la{\langle}

\def\leftitem#1{\item{\hbox to\parindent{\enspace#1\hfill}}}

\def\qed{\hfill\break\rightline{$\bbox$}}

\def\ra{\rangle}

\def\sg{\sigma}

\def\sg2{\sigma^2}

\def\__{_{\infty}}

\begin{document}

\title{An approach to chaos in some mixed $p$-spin models}
\author{Wei-Kuo Chen\footnote{Department of Mathematics, University of California at Irvine, email: weikuoc@uci.edu.}
\and 
Dmitry Panchenko \footnote{Department of Mathematics, Texas A$\&$M University, email: panchenk@math.tamu.edu. Partially supported by NSF grant.}}

\maketitle

\begin{abstract}
We consider the problems of chaos in disorder and temperature for coupled copies of the mixed $p$-spin models. Under certain assumptions on the parameters of the models we will first prove a weak form of chaos by showing that the overlap is concentrated around its Gibbs average depending on the disorder and then obtain several results toward strong chaos by providing control of the overlap between two systems in terms of their Parisi measures. 
\end{abstract}


\section{Introduction and main results.}

The phenomenon of chaos in disorder and temperature in spin glasses was discovered in \cite{FH86} and \cite{BM87} and has been studied extensively in the context of various models in the physics literature (e.g. see \cite{Rizzo09} for a recent review). In recent years, several mathematical results have also been obtained. An example of chaos in external field for the spherical Sherrington-Kirkpatrick model was given in \cite{Pan+Tal07}, chaos in disorder for mixed $p$-spin models with even $p\geq 2$ and without external field was considered in \cite{Chatt08}, \cite{Chatt09} (among many other results) and a more general situation in the presence of external field was handled in \cite{Chen11}. In this paper we will develop an approach to chaos in disorder and temperature for mixed $p$-spin models which is based on a novel application of the Ghirlanda-Guerra identities \cite{GG98} used here to derive a new family of identities in the setting of the coupled systems. At the moment, our approach only works under certain assumptions on the parameters of the models but these new examples are still welcome considering paucity of the results in this direction. Given $N\geq 1$, let us consider pure $p$-spin Hamiltonians $H_{N,p}(\vsi)$ for $p\geq 1$ indexed by $\vsi\in\Sigma_N = \{-1,+1\}^N$, 
\begin{equation}
H_{N,p}(\vsi)
=
\frac{1}{N^{(p-1)/2}}
\sum_{1\leq i_1,\ldots,i_p\leq N}g_{i_1,\ldots,i_p} \sigma_{i_1}\ldots\sigma_{i_p},
\label{HNp}
\end{equation}
where random variables $(g_{i_1,\ldots,i_p})$ are standard Gaussian independent for all $(i_1,\ldots,i_p)$ and $p\geq 1$. 
The covariance of this Gaussian process can be easily computed and is given by 
\begin{equation}
\frac{1}{N}\,\e H_{N,p}(\vsi^1) H_{N,p}(\vsi^2) = \bigl(R(\vsi^1,\vsi^2)\bigr)^p,
\label{CovHp}
\end{equation}
where quantity $R(\vsi^1,\vsi^2) = N^{-1}\sum_{i=1}^N \sigma_i^1\sigma_i^2$ is called the overlap of spin configurations $\vsi^1,\vsi^2\in\Sigma_N.$  Let us define a mixed $p$-spin Hamiltonian by a linear combination 
\begin{equation}
H_N(\vsi) = \sum_{p\geq 1} \beta_p H_{N,p}(\vsi)
\label{HN}
\end{equation}
with coefficients $(\beta_p)$ that decrease fast enough to ensure that the process is well defined, for example,  $\sum_{p\geq 1} 2^p \beta_p^2<\infty$. The Gibbs measure $G_N(\vsi)$ on $\Sigma_N$ is defined by
\begin{equation*}
G_N(\vsi) = \frac{\exp H_N(\vsi)}{Z_N},
\end{equation*}
where the normalizing factor $Z_N$ is called the partition function. The behavior of the Gibbs measure is intimately related to the computation of the free energy $N^{-1}\log Z_N$ in the thermodynamic limit and, as a result, has been studied extensively since the groundbreaking work of G. Parisi in \cite{Parisi79}, \cite{Parisi80}. In particular, various physical properties of the Gibbs measure, such as ultrametricity and lack of self-averaging, implied by the choice of the replica matrix in the Parisi ansatz were discovered by the physicists in the eighties (see \cite{MPV} for detailed account). The  chaos problem, or ``chaotic nature of the spin-glass phase" \cite{BM87}, arose from the discovery that, in some models, small changes in temperature or disorder may result in dramatic changes in the location of the ground state with the energy $\max_\vsi H_N(\vsi)$, as well as the overall energy landscape and the organization of the pure states of the Gibbs measure $G_N$. One very basic way to define such instability of the Gibbs measure is to sample a vector of spin configurations $\vsi$ from $G_N$ and a vector $\vrho$ from the measure $G_N'$ corresponding to the perturbed parameters and consider their overlap $R(\vsi,\vrho).$ The fact that this overlap behaves very differently than the overlap $R(\vsi^1,\vsi^2)$ of two replicas $\vsi^1, \vsi^2$ sampled from the same measure $G_N$ indicates that the set of configurations in $\Sigma_N$ on which the Gibbs measure concentrates (the location of pure states) is affected significantly by a small change of parameters. A typical statement that one is looking for in this case is that the overlap $R(\vsi,\vrho)$ is concentrated near zero when the model has symmetry or, more generally, near a constant when the symmetry is broken, for example, in the presence of external field. Indeed, this behavior is quite different from a typical ``lack of self-averaging" when the overlap between $\vsi^1$ and $\vsi^2$ can take many different values for any realization of the disorder in the low temperature phase. Moreover, even if we could show that the overlap $R(\vsi,\vrho)$ concentrates near its Gibbs' average which depends on the disorder, this would already indicate some form of chaos for exactly the same reasons. This is precisely what we will show in the case of perturbations of the disorder and for some perturbations of the inverse temperature parameters $(\beta_p)$. Furthermore, under additional assumptions on the sequence $(\beta_p)$ we will provide stronger control of the overlap in terms of the Parisi measures of the two systems. 

We will consider two systems with Gibbs' measures $G_N^1$ and $G_N^2$ corresponding to the Hamiltonians $H_N^1(\vsi)$ and $H_N^2(\vrho)$ as in (\ref{HN}) for $\vsi,\vrho\in\Sigma_N$ defined in terms of possibly different sequences of parameters $(\beta_p^1)$ and $(\beta_p^2)$ and, again, possibly different Gaussian disorders $(g^1_{i_1,\ldots,i_p})$ and $(g^2_{i_1,\ldots,i_p})$ for $p\geq 1$. However, we will assume that all pairs $(g^1_{i_1,\ldots,i_p}, g^2_{i_1,\ldots,i_p})$ are jointly Gaussian and independent for all $(i_1,\ldots,i_p)$ and $p\geq 1$. We will denote by $(\vsi^l,\vrho^l)_{l\geq 1}$ an i.i.d. sequence of replicas from the measure $G_N^1\times G_N^2$ and by $\la\cdot\ra$ the Gibbs average with respect to $(G_N^1\times G_N^2)^{\otimes\infty}.$

\smallskip
\noindent
\textbf{Weak forms of chaos.} For $j\in\{1,2\}$ let us denote
\begin{equation*}
\I_j^e=\bigl\{p \in 2\Natural : \beta_p^j\not = 0\bigr\},\,\,\I_j^o=\bigl\{p\in 2\Natural-1 : \beta_p^j\not = 0\bigr\}
\end{equation*}
and let $\I_j=\I_j^e\cup\I_j^o.$ When we talk about chaos in disorder we will assume that the following condition about their correlation is satisfied for at least one $p\geq 1$, 
\begin{equation}
p\in \I_1\cap\I_2 \,\mbox{ and }\,
{\rm corr} (g^1_{i_1,\ldots,i_p}, g^2_{i_1,\ldots,i_p} ) = t_p \in [0,1) 
\label{Corrp}
\end{equation}
for all $(i_1,\ldots,i_p)$. Our first result yields a weak form of chaos in disorder.
\begin{theorem}\label{ThWD} If (\ref{Corrp}) holds for some even $p\geq 2$ then
\begin{equation}
\lim_{N\to \infty}\e \bigl\la \bigl(|R(\vsi,\vrho)|-\la |R(\vsi,\vrho)|\ra \bigr)^2 \bigr \ra=0.
\label{Weven}
\end{equation}
If (\ref{Corrp}) holds for some odd $p\geq 1$ then
\begin{equation}
\lim_{N\to \infty}\e \bigl\la \bigl(R(\vsi,\vrho)-\la R(\vsi,\vrho)\ra \bigr)^2 \bigr \ra=0.
\label{Wodd}
\end{equation}
\end{theorem}\noindent
For example, for pure $3$-spin model, the overlap $R(\vsi,\vrho)$ is concentrated around its Gibbs average $\la R(\vsi,\vrho)\ra$ and for pure $2$-spin (SK) model, the absolute value of the overlap $|R(\vsi,\vrho)|$ is concentrated around its Gibbs average $\la |R(\vsi,\vrho)|\ra$. If $t_p=1$ in (\ref{Corrp}) for all $p\geq 1$, we can prove a weak form of chaos in temperature under certain assumptions on the sequences $(\beta^1_p)$ and $(\beta^2_p).$ Let us introduce a family of subsets of natural numbers,
\begin{equation}
\C_0 = \Bigl\{\I\subseteq \Natural : \mbox{linear span of $(x^p)_{p\in\I}$ is dense in $\bigl(C[0,1],\|\cdot\|_{\infty}\bigr)$} \Bigr\}.
\label{C0}
\end{equation}
Let us define the following conditions on the sequences $(\beta_p^1)$ and $(\beta_p^2)$:
\begin{enumerate}
\item[(${\rm C}_1^e$)] either $\I_2^e\setminus \I_1^e\not = \emptyset$ or there exist $\A\subseteq \I_1^e$ and  $p_0\in \I_1^e\setminus \A$ such that $\A\in\C_0$ and for some $\tau\in\Reals$ we have $\beta_p^2 = \tau \beta_p^1$ for $p\in\A$ and $\beta_{p_0}^2 \not = \tau \beta_{p_0}^1$,

\item[(${\rm C}_1^o$)] either $\I_2^o\setminus \I_1^o\not = \emptyset$ or there exist $\A\subseteq \I_1^o$ and  $p_0\in \I_1^o\setminus \A$ such that $\A\in\C_0$ and for some $\tau\in\Reals$ we have $\beta_p^2 = \tau \beta_p^1$ for $p\in\A$ and $\beta_{p_0}^2 \not = \tau \beta_{p_0}^1$,
\end{enumerate}
and let us define (${\rm C}_2^e$) and (${\rm C}_2^o$) in the same way by flipping indices $\{1, 2\}$. We will define conditions
\begin{equation}
({\rm C}^o)  = ({\rm C}^o_1) \wedge ({\rm C}^o_2),\,\,
({\rm C}^e) = (({\rm C}^e_1) \vee ({\rm C}^o_1))\wedge (({\rm C}^e_2) \vee ({\rm C}^o_2)).
\label{jointC}
\end{equation}
The role of (\ref{C0}) and condition $\A\in\C_0$ will be to ensure the validity of the extended Ghirlanda-Guerra identities from the identities for moments. The following weak form of chaos in temperature holds.
\begin{theorem}\label{ThWT} Condition {\rm(${\rm C}^e$)} implies (\ref{Weven}) and condition {\rm(${\rm C}^o$)} implies (\ref{Wodd}).
\end{theorem}\noindent
\textit{Example 1.} If $\I_1=\{3\}$ and $\I_2 = \{5\}$ then (\ref{Wodd}) holds. If $\I_1=\{2\}$ and $\I_2 = \{4\}$ or $\I_2 = \{3\}$ then (\ref{Weven}) holds.

\smallskip
\noindent
\textit{Example 2.} If $\I_1 = \I_2 = 2\Natural$,  $\beta_2^1 = \beta_2^2$ and $\tau\beta_p^1 = \beta_p^2$ for all even $p\geq 4$ and $\tau\not = 1$ then (\ref{Weven}) holds. 

\smallskip
\noindent
\textit{Example 3.} If $\I_1 = \{2\}$ and  $\I_2 = 2\Natural+2$, then (\ref{Weven}) holds.

\smallskip
\noindent
\textbf{Toward strong chaos.} To formulate the results that provide some strong control of the overlap $R(\vsi,\vrho)$ we need to recall some consequences of the validity of the Parisi formula for the free energy in mixed $p$-spin models, which was proved in \cite{Tal06} for even-spin models using the replica symmetry breaking interpolation idea from \cite{Guerra03} and in \cite{Pan11:1} in the general case using ultrametricity result from \cite{Pan11:4}. 
The first consequence that was found in \cite{Tal07} (see \cite{Pan11:2} or \cite{Tal11} for a simplified proof) states that the Parisi formula is differentiable in the inverse temperature parameters $\beta_p$ for all $p\geq 1$ which together with convexity implies that for all $p\in\I_1$,
\begin{equation}
\lim_{N\to\infty}\e \bigl\la R^p(\vsi^1,\vsi^2) \bigr\ra = \int_0^1 \! q^p \,d\mu_1(q),
\label{mu1}
\end{equation}
where $\mu_1$ is any probability measure on $[0,1]$ that achieves the minimum in the variational problem that defines the Parisi formula. Any such $\mu_1$ is called a Parisi measure of the system. Similarly, for all $p\in\I_2$,
\begin{equation}
\lim_{N\to\infty}\e \bigl\la R^p(\vrho^1,\vrho^2) \bigr\ra = \int_0^1 \!q^p \,d\mu_2(q)
\label{mu2}
\end{equation}
for any Parisi measure $\mu_2$ of the second system. Another consequence of the Parisi formula will be the strong form of the Ghirlanda-Guerra identities proved in \cite{Pan10:1} that will be used in the next section. In the situations that we  consider below the linear span of $(x^p)_{p\in\I_j}$ will be dense in $(C[0,1],\|\cdot\|_{\infty})$ for one or both $j=1,2$ in which case (\ref{mu1}), (\ref{mu2}) imply that the Parisi measure $\mu_j$ is unique. In this case let 
\begin{equation*}
c_j = \inf {\rm supp}\, \mu_j
\end{equation*}
be the smallest point in the support of $\mu_j$. The following result provides some control of the overlap and points toward strong chaos in disorder. 
\begin{theorem}\label{ThSD} If (\ref{Corrp}) holds for some $p\geq 1$ and $\I_j^e \in \C_0$ for $j=1$ or $j=2$ then
\begin{equation}
\lim_{N\to\infty}\e \bigl\la I\bigl(|R(\vsi,\vrho)| > \sqrt{c_j}\bigr) \bigr\ra=0.
\label{Seven}
\end{equation}
If (\ref{Corrp}) holds for some odd $p\geq 1$ and $\I_j^e \in \C_0$ for both $j=1$ and $j=2$ then
\begin{equation}
\lim_{N\to\infty}\e \bigl\la I\bigl(|R(\vsi,\vrho)| > \sqrt{c_1 c_2}\bigr) \bigr\ra=0.
\label{Sodd}
\end{equation}
\end{theorem}\noindent
In particular, if the Parisi measure $\mu_j$ of the system that satisfies $\I_j^e \in \C_0$ contains zero in its support then the overlap $R(\vsi,\vrho)$ concentrates around zero. Again, if $t_p=1$ in (\ref{Corrp}) for all $p\geq 1$, we have a similar result for chaos in temperature. 
\begin{theorem}\label{ThST} If $\I_j^e \in \C_0$ for $j=1$ or $j=2$ then condition {\rm(${\rm C}^e$)} implies (\ref{Seven}), and if $\I_j^e \in \C_0$ for both $j=1$ and $j=2$ then condition {\rm(${\rm C}^o$)} implies (\ref{Sodd}).
\end{theorem}\noindent
Finally, all our results also hold for the spherical mixed $p$-spin models when $\Sigma_N$ is the sphere of radius $\sqrt{N}$ with uniform measure, as long as $\I_1\cup \I_2\subseteq 2\Natural\cup\{1\}.$ This restriction is due to the fact that the Parisi formula for the spherical model has so far been proved only for such models in \cite{Tal061}.

\section{Ghirlanda-Guerra identities for coupled systems.}

In this section we will show how one can use the Ghirlanda-Guerra identities for each system in the form of the concentration of the Hamiltonian to obtain a new set of identities for the overlaps of the coupled system.  First of all, condition (\ref{Corrp}) means that the Gaussian pair $(g^1, g^2)$ is equal in distribution to $$(\sqrt{t_p}g + \sqrt{1-t_p}z^1,\sqrt{t_p} g + \sqrt{1-t_p}z^2)$$ for three independent standard Gaussian random variables $g,z^1,z^2$ and, therefore, the pair of processes $H_{N,p}^1(\vsi)$ and $H_{N,p}^2(\vrho)$ is equal in distribution to the pair
\begin{equation*}
\sqrt{t_p}H_{N,p}(\vsi) + \sqrt{1-t_p} Z_{N,p}^1(\vsi) \,\mbox{ and }\,
\sqrt{t_p} H_{N,p}(\vrho) + \sqrt{1-t_p} Z_{N,p}^2(\vrho),
\end{equation*}
where we denote by $H_{N,p}, Z_{N,p}^1$ and $Z_{N,p}^2$ three independent copies of (\ref{HNp}). Let us consider the quantities
\begin{align*}
\Gamma_{p}^1&=\e\Bigl\la \Bigl |\frac{H_{N,p}^1(\vsi^1)}{N}-\e \Bigl\la \frac{H_{N,p}^1(\vsi^1)}{N}\Bigr\ra \Bigr| \Bigr\ra,\\
\Gamma_{p}^2&=\e\Bigl\la \Bigl |\frac{H_{N,p}^2(\vrho^1)}{N}-\e \Bigl\la \frac{H_{N,p}^2(\vrho^1)}{N}\Bigr\ra \Bigr| \Bigr\ra,\\
\Delta_{p}^1&=\e\Bigl\la \Bigl |\frac{Z^2_{N,p}(\vsi^1)}{N}-\e \Bigl\la \frac{Z^2_{N,p}(\vsi^1)}{N}\Bigr\ra \Bigr| \Bigr\ra,\\
\Delta_{p}^2&=\e\Bigl\la \Bigl |\frac{Z^1_{N,p}(\vrho^1)}{N}-\e \Bigl\la \frac{Z^1_{N,p}(\vrho^1)}{N}\Bigr\ra \Bigr| \Bigr\ra.
\end{align*}
The Ghirlanda-Guerra identities \cite{GG98} in the form of the concentration of the Hamiltonians can be stated as follows. 
\begin{lemma}\label{Lem1add}
For all $p\geq 1,$ we have $\Delta_p^1, \Delta_p^2, \Gamma_p^1, \Gamma_p^2 \to 0.$
\end{lemma}\noindent
\textbf{Proof.}
Notice that in the definition of $\Delta_p^1$ we are averaging the Hamiltonian $Z^2_{N,p}$ from the second system over the first coordinate $\vsi^1$, which means that it is independent of the randomness in $\la\cdot\ra$. Therefore, if we denote by $\e'$ the expectation with respect to the randomness $Z_{N,p}^2$ then $\e\la Z_{N,p}^2(\s^1)\ra=\e \la \e' Z_{N,p}(\s^1)\ra=0$ and, using (\ref{CovHp}) and Jensen's inequality, 
\begin{equation*}
\e\Bigl\la\Bigl|\frac{Z^2_{N,p}(\vsi^1)}{N}\Bigr|\Bigr\ra\leq \e\Bigl\la\frac{\e' |Z^2_{N,p}(\vsi^1)|^2 }{N^2}\Bigr\ra^{1/2}\leq N^{-1/2}.
\end{equation*}
We conclude that $\Delta_p^1\to 0$ and, similarly, $\Delta_p^2\to 0.$ On the other hand, as we mentioned in the introduction, the validity of the Parisi formula for the free energy  and the argument in \cite{Pan10:1} (see also Chapter 12 in \cite{Tal11}) imply that $\Gamma_{p}^1\to 0$ and $\Gamma_p^2\to 0$ which is the usual formulation of the Ghirlanda-Guerra identities in the form of the concentration of the Hamiltonian. 
\qed

\noindent
Given replicas $(\vsi^l,\vrho^{l})_{l\geq 1}$ let us denote by
\begin{equation*}
R^1_{l,l'} = R(\vsi^l, \vsi^{l'}),\,\, R^2_{l,l'} = R(\vrho^l, \vrho^{l'}),\,\,R_{l,l'} = R(\vsi^l, \vrho^{l'})
\end{equation*}
the overlaps within each system and between the two systems. Notice that with these notations the cross overlap is not symmetric, $R_{l,l'}\not = R_{l',l}$. Given integer $n\geq 1,$ a function $\psi\in C\left[-1,1\right]$ and a bounded measurable function $f$ of the overlaps $(R_{l,l'}^1)_{l, l'\leq n},$ $(R_{l, l'}^2)_{ l,l'\leq n}$ and $(R_{l,l'})_{ l,l'\leq n}$ on $n$ replicas, we define
\begin{align}
\Phi_{1,n}(f,\psi)&=\e\la f\psi(R_{1,n+1}^1)\ra-\frac{1}{n}\e\la f\ra\hspace{0.3mm} \e\la\psi(R_{1,2}^1)\ra-\frac{1}{n}\sum_{l=2}^n\e\la f\psi(R_{1,l}^1)\ra,
\label{Phione}
\\
\Psi_{1,n}(f,\psi)&=\e\la f\psi(R_{1,n+1})\ra-\frac{1}{n}\sum_{l=1}^n\e\la f\psi(R_{1,l})\ra,
\label{Psione}
\\
\Phi_{2,n}(f,\psi)&=\e\la f\psi(R_{1,n+1}^2)\ra-\frac{1}{n}\e\la f\ra\hspace{0.3mm} \e\la\psi(R_{1,2}^2)\ra-\frac{1}{n}\sum_{l=2}^n\e\la f\psi(R_{1,l}^2)\ra,
\label{Phitwo}
\\
\Psi_{2,n}(f,\psi)&=\e\la f\psi(R_{n+1,1})\ra-\frac{1}{n}\sum_{l=1}^n\e\la f\psi(R_{l,1})\ra.
\label{Psitwo}
\end{align}
Throughout the paper we will use the notation 
\begin{equation*}
\psi_{p}(x)=x^p.
\end{equation*}
The following lemma contains a computation based on the Gaussian integration by parts analogous to the one for the original Ghirlanda-Guerra identities \cite{GG98} for one system.
\begin{lemma} \label{Lem1}
For all $p\geq 1$ we have,
\begin{align}
&\sup_{\|f\|_\infty\leq 1}\Bigl|\beta_p^2\sqrt{1-t_p}\Psi_{1,n}(f,\psi_p)\Bigl|\leq \frac{\Delta_{p}^1}{n},\label{Lem1:psi:C1}\\
&\sup_{\|f\|_\infty\leq 1}\Bigl|\beta_p^1\sqrt{1-t_p}\Psi_{2,n}(f,\psi_p)\Bigl|\leq \frac{\Delta_{p}^2}{n},\label{Lem1:psi:C2}\\
&\sup_{\|f\|_\infty\leq 1}\Bigl|\beta_{p}^1\Phi_{1,n}(f,\psi_p)+\beta_p^2 {t_p}\Psi_{1,n}(f,\psi_p)\Bigr|\leq \frac{\Gamma_{p}^1}{n}, \label{Lem1:phi:C1}\\
&\sup_{\|f\|_\infty\leq 1}\Bigl|\beta_{p}^2 \Phi_{2,n}(f,\psi_p)+\beta_p^1{t_p}\Psi_{2,n}(f,\psi_p)\Bigr|\leq \frac{\Gamma_{p}^2}{n}.
\label{Lem1:phi:C2}
\end{align}
\end{lemma}\noindent
In particular, Lemma \ref{Lem1add} implies that all the quantities on the left hand side go to zero and, under certain assumptions on the parameters of the models, this will imply that some or all quantities in (\ref{Phione}) - (\ref{Psitwo}) go to zero. Equations (\ref{Phione}) and (\ref{Phitwo}) will yield the familiar Ghirlanda-Guerra identities, only now the function $f$ may depend on the overlaps of the two systems. Furthermore, equations (\ref{Psione}) and (\ref{Psitwo}) will provide important additional information about how the two systems interact with each other. 

\smallskip
\noindent
\textbf{Proof.} We will only show $(\ref{Lem1:psi:C1})$ and $(\ref{Lem1:phi:C1})$ since  the proof of $(\ref{Lem1:psi:C2})$ and $(\ref{Lem1:phi:C2})$ is similar. As usual, we begin by writing that for $\|f\|_{\infty}\leq 1,$
\begin{equation}
\Bigl|\e\Bigl\la\frac{Z_{N,p}^2(\vsi^1)}{N}f\Bigr\ra-\e\Bigl\la\frac{Z_{N,p}^2(\vsi^1)}{N}\Bigr\ra \e\bigl\la f\bigr\ra\Bigr|
\leq \Delta_{p}^1
\label{Prop2:proof:D1}
\end{equation}
and
\begin{equation}
\Bigl|\e\Bigl\la\frac{H_{N,p}^1(\vsi^1)}{N} f\Bigr\ra-\e\Bigl\la\frac{H_{N,p}^1(\vsi^1)}{N}\Bigr\ra \e\bigl\la f\bigr\ra\Bigr|
\leq 
\Gamma_{p}^1.
\label{Prop2:proof:D2}
\end{equation}
Using (\ref{CovHp}) and Gaussian integration by parts we get
\begin{equation*}
\e\Bigl\la \frac{Z^2_{N,p}(\vsi^1)}{N} f\Bigr\ra
=\beta_p^2\sqrt{1-t_p}\Bigl(\sum_{l=1}^n\e\bigl\la (R_{1,l})^p f\bigr\ra-n\e\bigl\la (R_{1,n+1})^p f\bigr\ra\Bigr).
\end{equation*}
and since $\e\la Z^2_{N,p}\ra = 0,$ $(\ref{Prop2:proof:D1})$ implies $(\ref{Lem1:psi:C1}).$ Similarly, using Gaussian integration by parts, 
\begin{equation*}
\e\Bigl\la\frac{H_{N,p}^1(\vsi^1)}{N}\Bigr\ra=\beta_p^1\bigl(1-\e\bigl\la (R_{1,2}^1)^p\bigr\ra \bigr)
\end{equation*}
and
\begin{align*}
\e\Bigl\la\frac{H_{N,p}^1(\vsi^1)}{N} f\Bigr\ra
=\,\, &
\beta_p^1\Bigl(\sum_{l=1}^n \e\bigl\la (R_{1,l}^1)^pf\bigr\ra-n\e\bigl\la (R_{1,n+1}^1)^pf\bigr\ra\Bigr)
\\
&
+\, \beta_p^2 t_p\Bigl(\sum_{l=1}^n \e\bigl\la (R_{1,l})^pf\bigr\ra-n\e\bigl\la (R_{1,n+1})^pf\bigr\ra\Bigr).
\end{align*}
Therefore, $(\ref{Prop2:proof:D2})$ implies $(\ref{Lem1:phi:C1})$ and this completes the proof.
\qed

\noindent
We will use Lemmas \ref{Lem1add} and \ref{Lem1} in combination with the following result.

\begin{lemma}
\label{Lem2}
Let $j\in\{1,2\}.$ Suppose that
\begin{align}
&\lim_{N\to \infty}\sup_{\|f\|_\infty\leq 1}|\Psi_{j,n}(f,\psi)|=0
\label{Lem2:eq1}
\end{align}
holds with $\psi=\psi_{p}$ for some $p\geq 1.$ If $p\geq 2$ is even then (\ref{Lem2:eq1}) also holds for all even $\psi\in C[-1,1]$
and if $p\geq 1$ is odd then (\ref{Lem2:eq1}) holds for all $\psi\in C[-1,1]$.
\end{lemma}\noindent
\textbf{Proof.} It suffices to prove the results for $j=1.$ For all $l\geq 2$ (using symmetry),
\begin{align*}
\e\bigl\la\bigl((R_{1,1})^{p}-(R_{1,l})^{p}\bigr)^2\bigr\ra
=
2\e\bigl\la (R_{1,1})^{2p} \bigr\ra
-2\e\bigl \la (R_{1,1})^{p}(R_{1,2})^{p}\bigr\ra
=-2\Psi_{1,1}(f,\psi_{p})
\end{align*}
by definition of $\Psi_{1,n}$ in (\ref{Psione}) with $n=1$ and $f=(R_{1,1})^{p}$. If $p\geq 2$ is even then using that $|x-y|^{p }\leq |x^{p }-y^{p }|$ for all $x,y\geq 0$ we can write
\begin{align}
\e \bigl\la \bigl| |R_{1,1}|-|R_{1,l}| \bigr| \bigr\ra 
&\leq 
\bigl(\e \bigl\la \bigl||R_{1,1}|-|R_{1,l}|\bigr|^{2p } \bigr\ra\bigr)^{1/2p }
\label{Lem2:proof:eq}
\\
&\leq
\bigl(\e \bigl\la \bigl( (R_{1,1})^{p }-(R_{1,l})^{p } \bigr)^{2} \bigr\ra\bigr)^{1/2p }
=
\bigl(-2\Psi_{1,1}(f,\psi_{p})\bigr)^{1/2p}.
\nonumber
\end{align}
If $(\ref{Lem2:eq1})$ holds for $\psi=\psi_p$, this implies that $|R_{1,l}| \approx |R_{1,1}|$ for all $l\geq 2$ and, therefore, $(\ref{Lem2:eq1})$ holds for all even $\psi\in C[-1,1]$.  Whenever $(\ref{Lem2:eq1})$ holds for $\psi= \psi_{p }$ and odd $p \geq 1$ we use the same argument and the fact that $|x-y|^{p }\leq 2^{p-1}|x^{p }-y^{p }|$ for all $x,y\in\Reals$ to show that $R_{1,l} \approx R_{1,1}$ for all $l\geq 2$ and, therefore, $(\ref{Lem2:eq1})$ holds for all $\psi\in C[-1,1].$
\qed

\noindent
We are ready to state several consequences of Lemmas $\ref{Lem1add}$ - $\ref{Lem2}$ under additional assumptions on the parameters of the models that appear in our main results. First, we consider the condition (\ref{Corrp}) that is used to prove weak chaos in disorder.

\begin{proposition}\label{Prop1}
Suppose that (\ref{Corrp}) holds for some $p\geq 1.$ For $j\in \{1,2\},$ if $p$ is even then (\ref{Lem2:eq1}) holds for all even $\psi\in C[-1,1]$ and if $p$ is odd then (\ref{Lem2:eq1}) holds for all $\psi\in C[-1,1]$.
\end{proposition}\noindent
\textbf{Proof.}
Since under $(\ref{Corrp})$, $\beta_p^1,\beta_p^2 \not = 0$ and $t_p<1,$ equations (\ref{Lem1:psi:C1}), (\ref{Lem1:psi:C2}) and Lemma \ref{Lem1add} imply that (\ref{Lem2:eq1}) holds with $\psi = \psi_p$ for both $j\in\{1,2\}.$ The statement follows from Lemma $\ref{Lem2}.$
\qed

\noindent
One can prove a similar result under the conditions (\ref{jointC}) that appear in the results concerning chaos in temperature.

\begin{proposition}\label{Prop2}
Suppose that $t_p=1$ for all $p\geq 1.$ For $j\in\{1,2\}$, condition $({\rm C}^e)$ implies (\ref{Lem2:eq1}) for all even $\psi\in C[-1,1]$ and condition $({\rm C}^o)$ implies (\ref{Lem2:eq1}) for all $\psi\in C[-1,1]$.
\end{proposition}\noindent
\textbf{Proof.} The result will follow immediately from the definition of (${\rm C}^e$) and (${\rm C}^o$) in (\ref{jointC}) if we can show that 
\begin{enumerate}
\item[(i)] $({\rm C}_1^e)$ implies (\ref{Lem2:eq1}) for $j=1$ and even $\psi\in C[-1,1]$,  

\item[(ii)] $({\rm C}_{1}^o)$ implies (\ref{Lem2:eq1}) for $j=1$ and all $\psi\in C[-1,1]$,  

\item[(iii)] $({\rm C}_2^e)$ implies (\ref{Lem2:eq1}) for $j=2$ and even $\psi\in C[-1,1]$,  

\item[(iv)] $({\rm C}_{2}^o)$ implies (\ref{Lem2:eq1}) for $j=2$ and all $\psi\in C[-1,1]$. 
\end{enumerate}
We will only prove (i) since all other cases can be treated similarly. Let us show that $({\rm C}_1^e)$ implies 
\begin{equation}\label{Prop2:proof:eq1}
\lim_{N\to \infty}\sup_{\|f\|_\infty\leq 1}|\Psi_{1,n}(f,\psi_{p_0})|=0
\end{equation} 
for some even $p_0\geq 2$ from which (\ref{Lem2:eq1}) for $j=1$ and even $\psi\in C[-1,1]$ follows from Lemma $\ref{Lem2}$. First, if we suppose that $\I_{2}^e\setminus\I_1^e\neq\emptyset$ then there exists some even $p_0\geq 2$ such that $\beta_{p_0}^2\not= 0$ and $\beta_{p_0}^1=0,$ and (\ref{Prop2:proof:eq1}) immediately follows from $(\ref{Lem1:phi:C1})$. Next, suppose that there exist $\A\subseteq \I_1^e$ and  $p_0\in \I_1^e\setminus \A$ such that $\A\in\C_0$ and for some $\tau\in\Reals$ we have $\beta_p^2 = \tau \beta_p^1$ for $p\in\A$ and $\beta_{p_0}^2 \not = \tau \beta_{p_0}^1.$ Since $\beta_{p_0}^1\not = 0$, let $\beta_{p_0}^2/\beta_{p_0}^1=:\tau'\not = \tau$. Lemma \ref{Lem1add} and equation (\ref{Lem1:phi:C1}) imply that 
\begin{equation}\label{Prop3:proof:eq2}
\lim_{N\to \infty}\sup_{\|f\|_\infty\leq 1}|\Phi_{1,n}(f,\psi_{p_0})+\tau' \Psi_{1,n}(f,\psi_{p_0})|=0
\end{equation} 
and for $p\in\A$ (using that $\beta_p^2 = \tau \beta_p^1$ and $\beta_p^1\not = 0$),
\begin{equation*}
\lim_{N\to \infty}\sup_{\|f\|_\infty\leq 1}| \Phi_{1,n}(f,\psi_{p})+\tau \Psi_{1,n}(f,\psi_{p})|=0.
\end{equation*} 
Since $\A\in\C_0,$ we can approximate $\psi_{p_0}$ uniformly by $\psi_p$ for $p\in\A$ to obtain
\begin{equation}\label{Prop3:proof:eq3}
\lim_{N\to \infty}\sup_{\|f\|_\infty\leq 1}| \Phi_{1,n}(f,\psi_{p_0})+\tau \Psi_{1,n}(f,\psi_{p_0})|=0.
\end{equation} 
Since $\tau'\not = \tau$, (\ref{Prop3:proof:eq2}) and (\ref{Prop3:proof:eq3}) again imply (\ref{Prop2:proof:eq1}) and, thus, $({\rm C}_1^e)$ implies (\ref{Lem2:eq1}) for $j=1$ and even $\psi\in C[-1,1]$.
\qed

\noindent
Now that we obtained control of quantities $\Psi_{j,n}$, equations (\ref{Lem1:phi:C1}) and (\ref{Lem1:phi:C2}) can be used to control $\Phi_{j,n}$.

\begin{proposition}\label{Prop3}
Suppose that (\ref{Corrp}) holds for some $p\geq 1.$  For $j\in\{1,2\},$ if $\I_j^e \in \C_0$ then 
\begin{equation}\label{Prop3:psi}
\lim_{N\to \infty}\sup_{\|f\|_\infty\leq 1}|\Phi_{j,n}(f,\psi)|=0
\end{equation}
for all even $\psi\in C[-1,1]$.
\end{proposition}\noindent
\textbf{Proof.} Let us only consider the case $j=1$. By Proposition \ref{Prop1}, (\ref{Lem2:eq1}) holds for all even $\psi\in C[-1,1]$ and, therefore, equation (\ref{Lem1:phi:C1}) and Lemma \ref{Lem1add} imply that 
\begin{equation*}
\lim_{N\to \infty}\sup_{\|f\|_\infty\leq 1}|\Phi_{1,n}(f,\psi_p)|=0
\end{equation*}
for all $p\in \I_1^e.$ Since $\I_1^e\in \C_0$, we can approximate even $\psi\in C[-1,1]$ by polynomials with powers $p\in \I_1^e$ and $(\ref{Prop3:psi})$ follows for $j=1$.
\qed

\noindent
Exactly the same proof using Proposition \ref{Prop2} instead of Proposition \ref{Prop1} gives the following.

\begin{proposition}\label{Prop4}
Suppose that $t_p=1$ for all $p\geq  1.$
For $j\in\{1,2\},$ if $\I_j^e \in \C_0$ then either condition {\rm(${\rm C}^e$)} or {\rm(${\rm C}^o$)} implies (\ref{Prop3:psi}).
\end{proposition}


\section{Proof of the main results.}
As an immediate consequence of Propositions $\ref{Prop1}$ and $\ref{Prop2}$ we get Theorems $\ref{ThWD}$ and $\ref{ThWT}$. 

\smallskip
\noindent
\textbf{Proof of Theorems \ref{ThWD} and \ref{ThWT}.}
Suppose that either (\ref{Corrp}) holds for some even $p\geq 2$ or condition {\rm(${\rm C}^e$)} holds. By Propositions $\ref{Prop1}$ and $\ref{Prop2},$ (\ref{Lem2:eq1}) holds for all even $\psi\in C[-1,1]$ for both $j\in\{1,2\}$ and (\ref{Lem2:proof:eq}) implies
\begin{align*}
\lim_{N\to\infty}\e \bigl \la\bigl( |R_{1,1}|-|R_{1,2}| \bigr)^2 \bigr\ra=0.
\end{align*}
An argument similar to (\ref{Lem2:proof:eq}) also gives 
\begin{align*}
\lim_{N\to\infty}\e \bigl\la\bigl(|R_{2,2}|-|R_{1,2}|\bigr)^2 \bigr\ra=0.
\end{align*}
Equation $(\ref{Weven})$ follows by writing
\begin{align*}
\e \bigl\la\bigl(|R_{1,1}|- \la|R_{1,1}| \ra \bigr)^2 \bigr\ra
& \leq 
\e \bigl\la\bigl(|R_{1,1}|-|R_{2,2}|\bigr)^2 \bigr\ra
\\
& \leq 
2\e \bigl\la\bigl(|R_{1,1}|-|R_{1,2}|\bigr)^2 \bigr\ra
+
2\e \bigl\la\bigl(|R_{2,2}|-|R_{1,2}|\bigr)^2 \bigr\ra.
\end{align*}
If either (\ref{Corrp}) holds for some odd $p\geq 1$ or condition {\rm(${\rm C}^o$)} holds then, by Propositions $\ref{Prop1}$ and $\ref{Prop2},$ (\ref{Lem2:eq1}) holds for all $\psi\in C[-1,1]$ and a similar argument yields $(\ref{Wodd})$. 
\qed

\noindent
Let us denote by $\mu_N$ the distribution of the array of all overlaps
\begin{equation}
(R_{l,l'}^1)_{l, l'\geq 1}, (R_{l, l'}^2)_{ l,l'\geq 1} \,\mbox{ and }\, (R_{l,l'})_{ l,l'\geq 1}
\label{overlaps2}
\end{equation} 
under the annealed Gibbs measure $\e (G_N^1 \times G_N^2)^{\otimes \infty}$. By compactness, the sequence $(\mu_N)$ converges weakly over subsequences but, for simplicity of notation, we will assume that $\mu_N$ converges weakly to the limit $\mu$. We will still use the notations (\ref{overlaps2}) to denote the elements of the overlap array in the limit and, again, for simplicity of notations we will denote by $\e$ the expectation with respect to measure $\mu$. For example, whenever (\ref{Prop3:psi}) holds, the measure $\mu$ will satisfy the Ghirlanda-Guerra identities
\begin{equation}
\e f \psi(R_{1,n+1}^j)  = \frac{1}{n}\e f\, \e\psi(R_{1,2}^j) + \frac{1}{n}\sum_{l=2}^n\e f\psi(R_{1,l}^j)
\label{GGmu}
\end{equation}
for all bounded measurable functions $f$ of the overlaps on $n$ replicas and even $\psi\in C[-1,1].$ Consequently, (\ref{GGmu}) also holds for all even bounded measurable functions $\psi.$ Similarly, (\ref{Weven}) implies that $\mu$-almost surely $|R_{l,l'}| = |R_{1,1}|$ and (\ref{Wodd}) implies that $\mu$-almost surely $R_{l,l'} = R_{1,1}$ for all $l,l'\geq 1.$ Given $\mu$, let $\mu_1, \mu_2$ and $\mu_{1,2}$ denote the distributions of $|R_{1,2}^1|, |R_{1,2}^2|$ and $R_{1,1}$ under $\mu$ correspondingly (we will abuse the notations since, indeed, below the distributions of $|R_{1,2}^1|, |R_{1,2}^2|$ will coincide with the Parisi measures in (\ref{mu1}), (\ref{mu2})). Given measurable sets $A_1,A_2\subseteq [0,1]$ and $A\subseteq [-1,1]$ let us define the events
\begin{equation}
B_n =\Bigl\{R_{1,1}\in A,|R_{l,l'}^1|\in A_1 \mbox{ for } l \not = l'\leq n , |R_{l,l'}^2|\in A_2 \mbox{ for } l \not = l'\leq n \Bigr\}
\label{Lem3:def1}
\end{equation}
and
\begin{equation}
C_n =\Bigl\{R_{1,1}\in A,|R_{l,l'}^1|\in A_1 \mbox{ for } l \not = l'\leq n+1 , |R_{l,l'}^2|\in A_2 \mbox{ for }  l \not = l'\leq n \Bigr\}
\label{Lem3:def2}
\end{equation}
The following lemma will be crucial in the proof of Theorems \ref{ThSD} and \ref{ThST}. \begin{lemma}\label{Lem3}
If $\mu$ satisfies (\ref{GGmu}) for $j=1$ and $A_2=[0,1]$ then
\begin{equation}\label{Lem3:stat:eq1}
\mu(C_n)\geq {\mu}_1(A_1)^{n} {\mu}_{1,2}(A).
\end{equation}
If $\mu$ satisfies  (\ref{GGmu}) for $j=2$ and $A_1=[0,1]$ then
\begin{equation}\label{Lem3:stat:eq2}
\mu(B_{n})\geq {\mu}_2(A_2)^{n-1} {\mu}_{1,2}(A).
\end{equation}
If $\mu$ satisfies  (\ref{GGmu}) for both $j=1$ and $j=2$ then
\begin{align}
\begin{split}\label{Lem3:eq3}
\mu(B_n)\geq \left({\mu}_1(A_1) {\mu}_2(A_2)\right)^{n-1}{\mu}_{1,2}(A).
\end{split}
\end{align}
\end{lemma}\noindent
\textbf{Proof.}
Let us prove the following claim: If $\mu$ satisfies (\ref{GGmu}) for $j=1$ then
\begin{equation}\label{Lem3:eq1}
\mu(C_n)\geq {\mu}_1(A_1)\mu(B_n)
\end{equation}
and if $\mu$ satisfies (\ref{GGmu}) for $j=2$ then
\begin{equation}\label{Lem3:eq2}
\mu(B_{n+1})\geq {\mu}_2(A_2)\mu(C_n).
\end{equation}
First, we prove $(\ref{Lem3:eq1})$. We will use a computation similar to Lemma 1 in \cite{Pan11:3}. For any $n\geq 1$ we can write
\begin{align}\label{Lem3:proof:eq1}
I_{C_n}&\geq I_{B_n}-\sum_{ l \leq n}I_{B_{n}}I(|R_{l,n+1}^1|\notin A_1).
\end{align}
For all $1\leq l  \leq n,$ equation (\ref{GGmu}) for $j=1$ implies (using symmetry)
\begin{align*}
\e I_{B_{n}}I(|R_{l,n+1}^1|\notin A_1) 
\,&=\,
\frac{1}{n} {\mu}_1(A_1^c) \mu(B_{n})
+\frac{1}{n}\sum_{l '\neq l }^n \e I_{B_{n}} I(|R_{l,l'}^1|\notin A_1) 
\\
\,&=\,
\frac{1}{n}{\mu}_1(A_1^c)\mu(B_{n})
\end{align*}
and, therefore, $(\ref{Lem3:eq1})$ follows from $(\ref{Lem3:proof:eq1})$. In order to prove $(\ref{Lem3:eq2})$, let us start with
\begin{equation}
\label{Lem3:proof:eq2}
I_{B_{n+1}}\geq I_{C_n}-\sum_{l \leq n}I_{C_n}I(|R_{l,n+1}^2|\notin A_2).
\end{equation}
First of all, let us notice that  using the definition of the event $C_n$ and symmetry we can write for $l\leq n,$
$$
\e I_{C_n}I(|R_{l,n+1}^2|\notin A_2)  =\e I_{C_n}I(|R_{l,n+2}^2|\notin A_2) .
$$
Using (\ref{GGmu}) for the right hand side with $j=2$ and $n+1$ instead of $n$ (notice that $C_n$ depends on the first $n+1$ replicas),
\begin{align*}
\e I_{C_n}I(|R_{l,n+1}^2|\notin A_2) 
&=\frac{1}{n+1} {\mu}_2(A_2^c)\mu(C_n)+\frac{1}{n+1}\sum_{l '\neq l }^{n+1}
\e I_{C_n}I(|R_{l,l'}^2|\notin A_2)
\\
&=\frac{1}{n+1}{\mu}_2(A_2^c)\mu(C_n)+\frac{1}{n+1}\e I_{C_n}I(|R_{l,n+1}^2|\notin A_2).
\end{align*}
Therefore, for $1\leq l \leq n,$
$$
\e I_{C_n}I(|R_{l,n+1}^2|\notin A_2) 
=
\frac{1}{n}{\mu}_2(A_2^c)\mu(C_n)
$$
and $(\ref{Lem3:eq2})$ follows from $(\ref{Lem3:proof:eq2})$. Now suppose that (\ref{GGmu}) holds for $j=1$ and $A_2=[0,1]$. In this case, $C_n=B_{n+1}$ and $\mu(C_1)={\mu}_1(A_1){\mu}_{1,2}(A)$ using (\ref{GGmu}) with $n=1$.  
By induction, inequality (\ref{Lem3:eq1}) yields
$$
\mu(C_n)\geq {\mu}_1(A_1)^{n-1}\mu(C_1)= {\mu}_1(A_1)^n {\mu}_{1,2}(A)
$$ 
which proves $(\ref{Lem3:stat:eq1})$.  Now, suppose that (\ref{GGmu}) holds for $j=2$ and $A_1=[0,1].$ Then $C_n=B_n$ and by induction (\ref{Lem3:eq2}) implies (\ref{Lem3:stat:eq2}). Finally, suppose that (\ref{GGmu}) holds for both $j=1$ and $j=2$ and let us prove $(\ref{Lem3:eq3})$ by induction. First, it us easy to see that $\mu(B_1)={\mu}_{1,2}(A).$ Suppose that $(\ref{Lem3:eq3})$ holds for some $n\geq 1.$ Then using (\ref{Lem3:eq2}), (\ref{Lem3:eq1}) and induction hypothesis,
\begin{align*}
\mu(B_{n+1})
&\geq \,
{\mu}_2(A_2){\mu}(C_{n})
\geq 
{\mu}_1(A_1){\mu}_2(A_2)\mu(B_{n})
\\
&\geq \,
{\mu}_1(A_1){\mu}_2(A_2)\left({\mu}_1(A_1){\mu}_2(A_2)\right)^{n-1}{\mu}_{1,2}(A)
=
\left({\mu}_1(A_1){\mu}_2(A_2)\right)^n{\mu}_{1,2}(A).
\end{align*}
This completes the proof.
\qed


\noindent
\textbf{Proof of Theorems \ref{ThSD} and \ref{ThST}.}
Let us first prove the first part of the statements of Theorems \ref{ThSD} and \ref{ThST}. For certainly, let us assume that $\I_1^e\in\C_0.$ In this case (\ref{mu1}) implies that in the limit the distribution of $|R_{1,2}^1|$ coincides with the unique Parisi measure $\mu_1$. Suppose that either $(\ref{Corrp})$ holds for some $p\geq 1$ or, if not, condition $({\rm C}^e)$ holds. Then Propositions $\ref{Prop3}$ and $\ref{Prop4}$ imply that the identities (\ref{GGmu}) holds for $j=1$. Moreover, as we mentioned above, Theorems $\ref{ThWD}$ and $\ref{ThWT}$ imply that 
\begin{equation}\label{proof34:eq1}
\mu\bigl(|R_{1,1}|=|R_{l,l'}|,\forall l,l'\geq 1\bigr)=1.
\end{equation}
Let us show that the identities (\ref{GGmu}) for $j=1$ together with $(\ref{proof34:eq1})$ imply $(\ref{Seven})$, that is, 
\begin{equation}\label{Thm3:proof:eq1}
\lim_{N\rightarrow\infty}\mathbb{E}\left<I(|R_{1,1}^1|>\sqrt{c_1})\right>=0.
\end{equation}
Suppose that $(\ref{Thm3:proof:eq1})$ is not true. Then there exists some $c>\sqrt{c_1}$ such that $${\mu}_{1,2}(\left[-1,-c\right)\cup\left(c,1\right])>0.$$ Since $c_1$ is the smallest value of the support of $\mu_{1},$ there exists some $c_0$ with $c_1<c_0<c^2$ such that $\mu_1(\left[0,c_0\right))>0.$ Set $A=\left[-1,-c\right)\cup\left(c,1\right]$, $A_1=\left[0,c_0\right),$ and $A_2=\left[0,1\right].$ Recall the definition of $C_n$ from $(\ref{Lem3:def2}).$  Using $(\ref{Lem3:stat:eq1}),$ we know that $\mu(C_n)\geq ({\mu}_1(A_1))^n{\mu}_{1,2}(A)>0$ for each $n\geq 1.$ Let us consider the event
$$
\hat{C}_n =\bigl\{R_{l,1}\in A \mbox{ for } l\leq n+1,|R_{l,l'}^1|\in A_1 \mbox{ for } l \not = l'\leq n+1 \bigr\}.
$$
By $(\ref{proof34:eq1})$, $\mu(\hat{C}_n)=\mu(C_n)$ and, since $\hat{C}_n$ is an open subset on the space of overlaps, 
\begin{equation*}
\liminf_{N\rightarrow\infty}\mu_N(\hat{C}_n)\geq \mu(\hat{C}_n)=\mu(C_n)>0.
\end{equation*}
This means that for any $n\geq 2,$ for large enough $N$, we can find $\vsi^1,\vsi^2,\ldots,\vsi^{n}\in \Sigma_{N}$ and $\vrho^1\in \Sigma_N$ such that $|R(\vsi^l,\vrho^1)|>c$ for $l\leq n$ and $|R(\vsi^l,\vsi^{l'})|< c_0$ for $l \not = l'\leq n$. Let us choose $a_1,\ldots,a_{n}\in\left\{-1,1\right\}$ such that $a_l R(\vsi^l,\vrho^1)=|R(\vsi^l,\vrho^1)|$ for $l\leq n.$
Then
\begin{align*}
N^{-1}
\left|(a_1\vsi^1+a_2\vsi^2+\cdots+a_{n}\vsi^{n},\vrho^1)\right|
&=
\sum_{1\leq l\leq n} \bigl|R(\vsi^l,\vrho^1)\bigr|\geq nc
\end{align*}
and
\begin{align*}
N^{-1}
\bigl\|a_1\vsi^1+a_2\vsi^2+\cdots+a_{n}\vsi^{n}\bigr\|^2
&=
\sum_{1\leq l,l'\leq n}a_{l} a_{l'}R(\vsi^{l},\vsi^{l'})
\leq n+ (n^2-n)c_0.
\end{align*}
Using the Cauchy-Schwarz inequality, we obtain that
\begin{align*}
n^2c^2&\leq N^{-2} \bigl|(a_1\vsi^1+a_2\vsi^2+\cdots+a_{n}\vsi^{n},\vrho^1)\bigr|^2\\
&\leq N^{-2} \bigl\|a_1\vsi^1+a_2\vsi^2+\cdots+a_{n}\vsi^{n}\bigr\|^2 \|\vrho^1\|^2
\leq n+(n^2-n)c_0.
\end{align*}
If we divide both sides by $n^2$ and let $n\to\infty$ we get $c^2\leq c_0$ which contradicts the choice of $c_0.$ This completes the proof of (\ref{Seven}). Next, we prove (\ref{Sodd}) assuming that $\I_j^e\in \C_0$ for both $j=1$ and $j=2.$
In this case, the Parisi measures $\mu_1$ and $\mu_2$ are again the limiting distributions of $|R_{1,2}^1|$ and $|R_{1,2}^2|$, respectively. Suppose that either (\ref{Corrp}) holds for some odd $p\geq 1$ or $({\rm C}^o)$ holds. By Propositions
\ref{Prop3} and \ref{Prop4}, the identities (\ref{GGmu}) are satisfied for both $j=\{1, 2\}$ and, by Theorems $\ref{ThWD}$ and $\ref{ThWT},$ 
\begin{equation}\label{proof34:eq2}
\mu\bigl(R_{1,1}=R_{l,l'},\forall l,l'\geq 1\bigr)=1.
\end{equation}
We prove $(\ref{Sodd})$ by contradiction. Assume that there exists some $c>\sqrt{c_1c_2}$ such that 
$$
{\mu}_{1,2}\left(\left[-1,-c\right)\cup\left(c,1\right]\right)>0.
$$
Let us discuss the case ${\mu}_{1,2}(\left(c,1\right])>0$ first.
Choose $d_1$ and $d_2$ satisfying $c_1<d_1<1$, $c_2<d_2<1$, and $\sqrt{d_1d_2}<c$. 
If we define $A_1=[0,d_1)$, $A_2=[0,d_2)$  and $A=(c,1]$ then $\mu_1(A_1)>0$ and $\mu_2(A_2)>0.$ 
If we recall the event $B_n$ in (\ref{Lem3:def1}), (\ref{Lem3:eq3}) implies that  $\mu(B_n) >0$. If we consider the event
\begin{align*}
\hat{B}_{n}
=
\bigl\{
R_{l,l'}\in A \mbox{ for } l,l'\leq n, |R_{l,l'}^1|\in A_1,
|R_{l,l'}^2|\in A_2 \mbox{ for } l\neq l'\leq n
\bigr\}
\end{align*}
then by $(\ref{proof34:eq2})$, $\mu(\hat{B}_n)=\mu(B_n),$ and since $\hat{B}_n$ is an open subset on the space of overlaps, 
\begin{equation*}
\liminf_{N\rightarrow\infty}\mu_N(\hat{B}_n) \geq \mu(\hat{B}_n)=\mu(B_n)>0.
\end{equation*}
This implies that for any $n\geq 2$, if $N$ is sufficiently large, we can find $\vsi^1,\ldots,\vsi^n\in\Sigma_N$ and
$\vrho^1,\ldots,\vrho^n\in \Sigma_{N}$ such that $R(\vsi^{l},\vrho^{l'})\in A$ for $l,l'\leq n,$ $|R(\vsi^{l},\vsi^{l'})|\in A_1$ for $l\neq l'\leq n$ and $|R(\vrho^{l},\vrho^{l'})|\in A_2$ for $l\neq l'\leq n.$ Therefore,
\begin{align*}
N^{-1}\bigl\|\vsi^1+\vsi^2+\cdots+\vsi^{n}\bigr\|^2&=\sum_{l,l'\leq n}R(\vsi^{l},\vsi^{l'})\leq  n+(n^2-n)d_1,\\
N^{-1}\bigl\|\vrho^1+\vrho^2+\cdots+\vrho^{n}\bigr\|^2&=\sum_{l,l'\leq n}R(\vrho^{l},\vrho^{l'})\leq n+(n^2-n)d_2
\end{align*}
and
\begin{align*}
N^{-1}\bigl|\bigl(\vsi^1+\vsi^2+\cdots+\vsi^{n},\vrho^1+\vrho^2+\cdots+\vrho^{n}\bigr)\bigr|
=\Bigl|\sum_{l,l'\leq n}R(\vsi^{l},\vrho^{l'})\Bigr|\geq n^2c.
\end{align*}
Using the Cauchy-Schwarz inequality as above, 
\begin{align*}
n^4c^2\leq (n+(n^2-n)d_1) (n+(n^2-n)d_2).
\end{align*}
Since this is true for every $n,$ dividing both sides by $n^2$ and passing to the limit, it implies $c\leq \sqrt{d_1d_2}$
which contradicts the choice of $d_1$ and $d_2.$ This completes the proof in the case ${\mu}_{1,2}(\left(c,1\right])>0.$ 
One can check that the same argument yields the result when ${\mu}_{1,2}(\left[-1,-c\right))>0$ and this finishes the proof. 
\qed


\end{document}